# The TVBG-SEIR spline model for analysis of COVID-19 spread, and a Tool for prediction scenarios

Ognyan Kounchev,

Institute of Mathematics and Informatics, Bulgarian Academy of Sciences

Georgi Simeonov,

Institute of Mathematics and Informatics, Bulgarian Academy of Sciences

and

Zhana Kuncheva, C4X Discovery Ltd., UK

**ABSTRACT.**

The coronavirus infection (COVID-19) has started in December 2019 in Wuhan, China, and it has quickly spread out to almost all other countries in the world. Mathematical models are traditionally used to analyze the long-term global evolution of epidemics, to determine the potential and severity of an outbreak, and to provide critical information for identifying the type of disease interventions and intensity. One of the widely used mathematical models of long-term spreading of epidemics are the so-called deterministic compartmental models (SIR/SEIR type models) which reflect adequately the clinical progression of the virus spreading, epidemiological status of the individuals, and the intervention measures. One of the main purposes of applying such models is to assess how the expensive restriction measures imposed by the authorities (home and social isolation/quarantine, travel restrictions, etc.) can effectively reduce the control reproduction number of the disease and its transmission risk. In particular, it is essential to assess how the expensive, resource-intensive measures implemented by the authorities can contribute to the prevention and control of the COVID-19 infection, and how long they should be maintained.

However the classical SIR/SEIR models have been primarily studied in what may be called stationary case, where the main parameters, the Transmission rate β (reflecting the virus spread by infected individuals) and the Removed rate γ (reflecting the hospitalization/isolation measures) remain constant during the whole period of interest. This does not reflect in a proper way their extremely dynamical behavior during the COVID-19 (and similar) epidemic, resulting from the imposition of intensive restriction measures by the authorities.

Hence, it is important to extend the classical SIR/SEIR models by creating new ansatzes for the dynamics of the transmission rates β(t) (which we will call further just **Beta**) and removed rates γ(t) (which we will call further just **Gamma**). The main purpose of the present research is to introduce a spline-based SEIR model with Time-varying Beta and Gamma parameters, or abbreviated TVBG-SEIR model, which is used to estimate the practical implications of the public health interventions and measures. We have designed a Tool based on the TVBG-SEIR model, which may be used as a Decision Support Tool to assist the health decision- and policy-makers in creating predictive scenarios. It may be used to assess the impact of previous public health interventions, and to plan quantitatively and qualitatively the introduction of future containment measures for achieving the necessary objectives.

## Aim and Methods summary

To estimate the dynamics of transmission β and removed rates γ during a controlled COVID-19 outbreak we have developed a mathematical model with time-varying β(t) and γ(t) rates, TVBG-SEIR, which was simultaneously fitted to two sets of data: the daily infected cases, and the daily removed cases.

We use deterministic spline Ansatz: the transmission rates β(t) and the removal rates γ(t) are modeled by *splines* with two nodes - Node1, Node2 (the same nodes for both) - within the time interval of interest – from Start Date until Today date. This Ansatz allows to properly model the dynamics due to the introduction of containment measures by the authorities. The purpose of fitting of the TVBG-SEIR model is to identify the nodes of the splines and the three values of β(t) and γ(t) on Node1, Node2, and Today date. It is assumed that β(t) and γ(t) are constant in the time interval [Startdate, Node1], and β(t) is monotone decreasing while γ(t) is monotone increasing function.

A Tool was designed for the visualization of the results of the fitted model (the daily infected cases), and for creating prediction scenarios for the daily infected cases during the next two months, by controlling the future values of Beta and Gamma. It is available at the link:

http://213.191.194.141:8888/notebooks/TVBG-SEIR-Spline-model_v3.ipynb?token=b5d97bfbd7dd062e47ee7ab51837e470a8c226743a4667ee

The plan of the paper is as follows: In Section 1 we recall the deterministic SEIR model and introduce some notions and notations. In Section 2 we introduce the discretization of the SEIR model which is used in the algorithms. In Section 3 we introduce and provide all technical details of the TVBG-SEIR spline model. In Section 4 we provide an application of the TVBG-SEIR model to Bulgarian data, which are used to illustrate the work of the TOOL for prediction scenarios. In Section 5 we do the same for Italian data, and in Section 6 for German data. In Section 7 we describe the technical details of the Tool for prediction scenarios. In Section 8 we provide some recent references about models with time-varying transmission rates and their calibration.

## Data Acknowledgement

We acknowledge the following data sources about COVID-19:

HDX Humanitarian Data Exchange V1.39.3
- time_series_covid19_confirmed_global.csv
- time_series_covid19_deaths_global.csv
- time_series_covid19_recovered_global.csv

https://data.humdata.org/dataset/novel-coronavirus-2019-ncov-cases

Kaggle SRK Sudalairaj Kumar
https://www.kaggle.com/sudalairajkumar/novel-corona-virus-2019-dataset#covid_19_data.csv
  • covid_19_data.csv

# 0. Introduction

In the case of the usual seasonal flu the main parameters of the spread of the viruses are the transmission rate $\beta$ which reflects the power of the transmission of the virus from infected people to susceptibles, the recovery rate $\gamma$ which is reciprocle to the recovery period (which is the sum of recovery to health + isolated sick people + mortality due to the sickness ), and the parameter $\sigma$ which is the reciprocle to the incubation period. Due to the long incubation period and large number of asymptomatic or mild-symptomatic cases, COVID-19 has proved to be very insidious and requires intensive emergency measures from the authorities to reduce the transmission rate $\beta$ and to increase the recovery rate $\gamma$. For comparison, in the case of the seasonal flu no intensive measures are necessary to be undertaken by the authorities.

The authorities have introduced very strong restrictive measures which have essentially influenced the dynamics of the parameters $\beta$ and $\gamma$. For the majority of the states these measures have been introduced not only in one step but most often at least in two steps. It is very dependent on every society how fast will these measures be implemented in life. The measures are, for example, closing schools, pubs, restaurants, social meetings, wearing masks, etc. All they restrict the contacts among the people, and thus essentially change the dynamics of the $\beta$ and $\gamma$ rates. It is important to assess how these expensive and resource intensive measures implemented by the authorities can contribute to the prevention and control of the COVID-19 infection, and how long they should be maintained, [12], [13].

In order to meet the challenge of Controlled spread of the COVID-19 (and similar) epidemics, one needs to develop new mathematical models which describe better the reality. Some of the most widely used models in Epidemiology are the deterministic SEIR model ([4]). Based on it, in the present research we propose a new model TVBG-SEIR which incorporates a specific spline model for the time-varying transmission and removal rates.

# 1. The classical deterministic SEIR model, Notions and Notations

Compartmental models are a framework used to model in an adequate way the dynamics of infectious disease (see the Wiki article). The population is divided into compartments, with the assumption that every individual in the same compartment has the same characteristics. This framework has been developed for the first time in the paper of Kermack and McKendrick in 1927 [2]. One may use a deterministic approach using a system of ODEs or a stochastic approach which is more complicated. The deterministic approach is what we follow and has two main representatives - the SIR and the SEIR models.

For a detailed and excellent introduction to the compartmental SIR/SEIR models we refer to the monograph [4]. We provide a short description of the deterministic SEIR model which will be the main approach in our research.

The classical SEIR model is based on the consideration of four compartments, $C_S$, $C_E$, $C_I$ and $C_R$ which are described as follows:

1. Compartment $C_S$ : its size is $Sr(t)$ – the number of "susceptible" people at time $t$. Usually at the start $S(0) = N$ is the whole population of the country under consideration. It is supposed that nobody has automatic immunity against the virus, i.e. everybody is susceptible.
2. Compartment $C_E$ : its size is $Er(t)$ – the number of "exposed" people at time $t$ - these are the people who are "virus carriers" but are not "virus spreaders"; the virus is in a latent form, and usually they do not show symptoms of sickness. For different viruses the incubation (latent) period is very different – for the coronavirus it was recently statistically estimated that the average incubation period is 11.5 days [3]. Not everybody in $C_E$ may become "virus spreader", i.e. move to the next compartment $C_I$ . Practically, the compartment $C_E$ does not enter the official statistics, and it is practically not an observable but is very important for a more adequate modeling of the dynamics of the virus spread. This compartment is missing in the simpler SIR model.
3. Compartment $C_I$ : its size is $Ir(t)$ – the number of infectious cases at time $t$ - these are the people who are "virus spreaders", majority of them show some symptoms, although they may not show any symptoms (asymptomatic). It is important to understand in the modeling that many people who are diagnosed positively are almost immediately hospitalized or quarantined, hence they go to compartment $C_R$, but they have stayed in $C_I$ only until they have been diagnosed (and these are the official data which we obtain – $Idata(t)$ ).
4. Compartment $C_R$ : its size is $Rr(t)$ – the number of recovered or deceased (or immune) individuals, which are all called "removed". Normally they come from compartment $C_I$ after becoming healthy and no more virus spreaders. Officially these data are provided in a cumulative way.

However, what data do we have at our disposal to fit the model?

We do not have the "reality data" $Sr(t), Er(t), Ir(t), Rr(t)$. We have the official data $Idata(t)$ which are the daily "new infected cases" with COVID-19, and these are normally people with serious symptomatic. These are the cases which have been tested and registered officially at the hospitals. The majority of them are almost immediately hospitalized or quarantined, hence, they are almost immediately moved from compartment $C_I$ to compartment $C_R$. However it is well known that for seasonal flu (and it is considered to be similar for COVID-19) the size of $C_I$ is much bigger than that indicated by the official data $Idata(t)$, and we have the inequality

$$Idata(t) \leq Ir(t).$$

We have also the officially announced data $Rdata(t)$ which contain the cumulative number of recovered cases, and the $Ddata(t)$ which is the cumulative number of fatalities. Although there is a lot of discussion about the quality of these data, it is approximately true that

$$Rdata(t) + Ddata(t) \leq Rr(t)$$

A main point of the modeling paradigm for COVID-19 (and similar virus infections) is that, for a certain segment of the society (in this case, the younger people), the infection symptoms do not differ essentially from a seasonal flu, hence the number of unreported cases (those which are in compartment $C_I$ but not in $Idata(t)$ for every time $t$) may be much bigger, thus in the above inequalities more appropriate is to use the symbol "$\ll$" , which denotes "much less". In the case of the seasonal flu it may be even 100 times less.

The main point of developing the compartmental deterministic SEIR model is to provide some tractable approximations $S(t), E(t), I(t), R(t)$ to the above time series of the "reality data" $Sr(t)$, $Er(t)$, $Ir(t)$, $Rr(t)$. The most widely used is the model based on a system of Ordinary Differential Equations with variables $S(t)$, $E(t)$, $I(t)$, $R(t)$ which is given as follows:

$$S'(t) = -\beta(t)S(t)I(t)/N \tag{1a}$$

$$E'(t) = \beta(t)S(t)I(t)/N - \sigma E(t) \tag{1b}$$

$$I'(t) = \sigma E(t) - \gamma(t)I(t) \tag{1c}$$

$$R'(t) = \gamma(t)I(t) \tag{1d}$$

Let us explain the notations:

1. Here the term $\beta(t)I(t)/N$ expresses the rate at which new individuals (as a proportion of the total population size) are infected by the already infectious $I(t)$ individuals. Here and further $\beta(t)$ is called Transmission rate of the infection, which we call further simply Beta.

2. The coefficient $\gamma(t)$ is the Removal or recovery rate; it is determined by the reciprocle of the infectious period, after which either the person is recovered (and no more infectious) or dead (again, no more infectious). Here and further $\gamma(t)$ is called Removal rate, and we call it simply Gamma.

3. The coefficient $\sigma$ is the latent rate, or the rate of "becoming symptomatic" (where $1/\sigma$ is the average of the incubation period). In the present paper we use the constant value

$$\sigma = 1/5.2$$

which represents a reasonable approximation, as the recent research shows, see [3].

The usual applications of the SEIR model are with constant rates $\beta(t)$ and $\gamma(t)$. One assumes that the initial values $S(0)$, $E(0)$, $I(0)$, and $R(0)$ are given and the system is solved for the times $t \geq 0$, where $t$ is an integer. It is assumed that the following "conservation" equation holds

$$N = S(t) + E(t) + I(t) + R(t) \tag{2}$$

where $N$ is the total population in the country $XX$. Obviously, after introducing equation (2), the fourth equation in (1) becomes redundant.

## 2. Discretization of the SEIR model

In practice one uses a discretization of the continuous SEIR model. The following discretization of the SEIR model is very intuitive, and is in fact derived from the Euler method for approximate solution of the initial value problem (1):

$$S_{n+1} = S_n - \beta_n S_n I_n / N \tag{3a}$$

$$E_{n+1} = E_n + \beta_n S_n I_n / N - \sigma E_n \tag{3b}$$

$$I_{n+1} = I_n + \sigma E_n - \gamma_n I_n \tag{3c}$$

$$R_{n+1} = R_n + \gamma_n I_n \tag{3d}$$

Here $S_n$, $E_n$, $I_n$ and $R_n$ are respectively the values of $S(t)$, $E(t)$, $I(t)$ and $R(t)$ on the day $t = n$, and the initial values for day $n = 0$ are $S_0$, $E_0$, $I_0$ and $R_0$. The above system is iteratively solved for integers $n \geq 0$. We assume that the size $N$ of the population remains unchanged (hence no usual birth and mortality are taken into account). Hence, the total sum of the above is assumed to satisfy

$$N = S_n + E_n + I_n + R_n \tag{3e}$$

which makes the fourth equation in (3) redundant.

It is well known that the above Euler method for approximating the solution of (1) is less accurate than the Runge-Kutta which is widely used, see e.g. [11].

**Remark 3** One has to note that the continuous model (1a)-(1d) and the above discrete approximation (3a)-(3d) have essential differences in the long-term behavior which has been the subject of much research. It is important to note that the qualitative properties of the solution to the differential equation and of the discrete equation differ essentially - the continuous case is simpler as usual.

## 3. The TVBG-SEIR model

The SIR/SEIR models have proved to be very efficient in situations where the main parameters $\beta$ and $\gamma$ are constants, in natural conditions, where no special control by the authorities is exercised, i.e. no intervention (containment) measures are undertaken to change the transmission and the removal rates in the course of the epidemics. This is very often the case with the seasonal flu where the medical authorities do not undertake actively special measures to restrict the social behaviour of the citizens. However due to the specific of the COVID-19 the situation has become more dramatic and it has required the interference of the governments in order to avoid the overloading of the National Health systems. The authorities have introduced very strong restrictive measures which have essentially influenced the dynamics of the parameters $\beta$ and $\gamma$. For the majority of the states these measures have been introduced not only in one step but most often at least in two steps.

In view of the above it makes sense to seek for Mathematical models which try to model as best as possible the dynamical change of the parameters $\beta$ and $\gamma$. We have decided for spline structure with two important breakpoints Node1, Node2 which reflect the control exercised by the authorities in the form of restriction measures. Also, it is natural to assume that inbetween the dates the control measures change the parameters $\beta(t)$ and $\gamma(t)$ in a monotone way, i.e. $\beta(t)$ is decreasing whereas $\gamma(t)$ is increasing.

Technical description of the TVBG-SEIR spline model:

1. We denote the Start date by **T1**; this corresponds to a date when the first cases of COVID-19 are announced, eventually we may choose **T1** to be a date when the steeper growth of the epidemic starts. We denote by **T4** the end date (usually chosen to be Today).
2. We choose two interior nodes for the *interpolation splines* modeling the coefficients Beta and Gamma: Node1 = **T2** and Node2 = **T3**. This corresponds to two steps of the introduction of Restrictive Measures imposed by the authorities of the country **XX**. Normally, the date **T2** may be the First restrictive measures date, or a date close to it, and **T3** may be the Second restrictive measures date, or a date close to it.

3. The model is supposed to reflect the natural expectation that once there are official restrictions they will implicate an essential change in the Transmission and Removed rates although not immediately. We assume that the Beta rate *β(t)* is monotone decreasing with the time, which corresponds to the natural expectation that the more restrictive the measures the smaller the Transmission rate. Respectively, the Gamma rate γ(t) is assumed to be monotone increasing, meeting the expectation that the stronger the measures, the bigger the removal rate.

4. We assume that *β(t)* and *γ(t)* are constant between the Start date **T1** and the first node **T2**, i.e. *β(T1) = β(T2)* and *γ(T1) = γ(T2)*. This corresponds to the "still" life of the society (without containment measures) when the Transmission and the Removal rates are nearly a constant.

5. To be more precise, the splines which we consider are not the usual polynomial, but the so-called Exponential splines depending on a parameter in the exponent, which makes a fast decay to the next target value of the Beta rates; respectively this makes fast increasing to the target value of the Gamma rates. This corresponds to the expectation that the speed by which the society switches from one level of the restrictive measures to another is relatively fast, and it is reflected by the size of the exponent we decide to choose. On the following Figures we see examples of the dynamics of *β(t)* and *γ(t)* rates:

Figure 1  This is the Beta coefficient for a Model curve

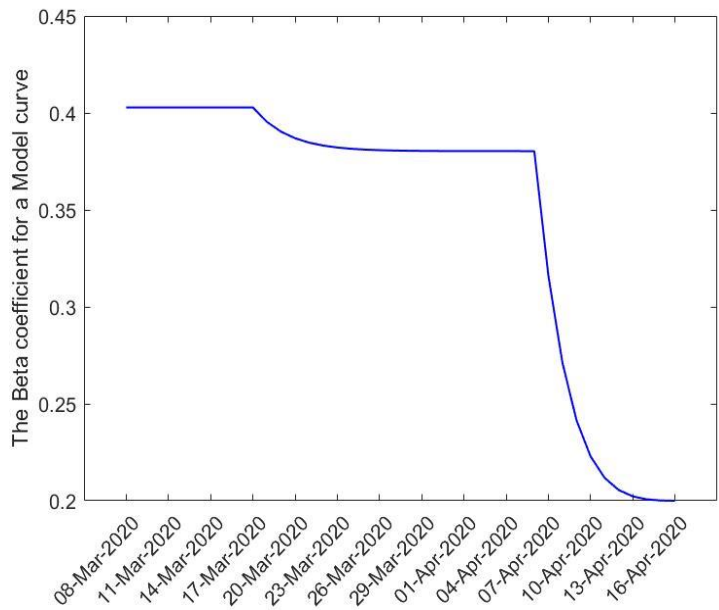

Figure 2  This is the Gamma coefficient for a Model curve

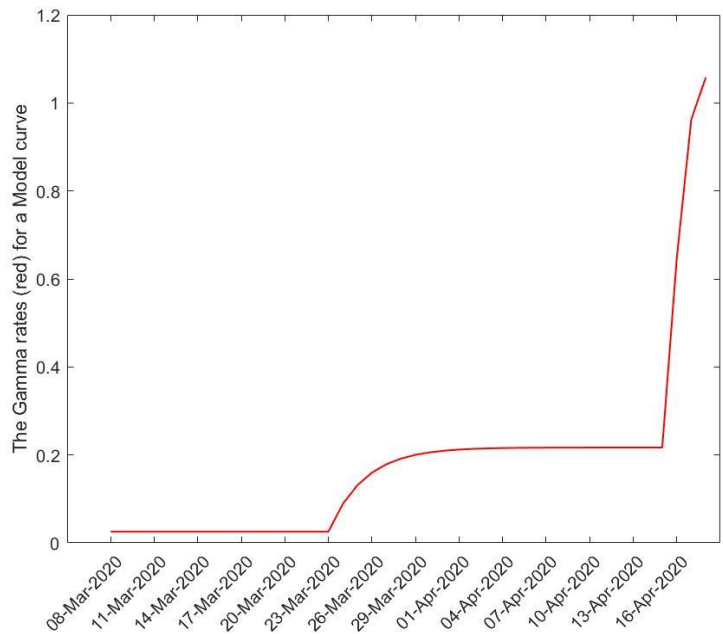

6. An important property of the TVBG-SEIR model is that due to the above spline model for the Beta and Gamma parameters, where there is a fast transition to the next target value, a classical SEIR model with constant β(t) and γ(t) holds during larger sub-intervals. In particular, this permits to provide a reliable estimate of the Basic Reproduction Number (Ratio).
7. The Reproduction number (ratio) is a key variable for all models of epidemics, see [4], [9], [14], [15]. Following [15] (formula (2.4)), for the case of the SEIR models with constant rates β(t) and γ(t), the Reproduction number is given by the formula

$$R_0 = \frac{\beta}{\gamma}$$

where we have assumed that the natural birth and mortality rates are small and also equal. Due to the above remark, we may extrapolate the above formula for all time points of interest by putting:

$$R_0 = \frac{\beta(t)}{\gamma(t)}$$

On the Figure below we provide the Reproduction Number (Ratio), obtained by the last formula, for some specific TVBG-SEIR model:

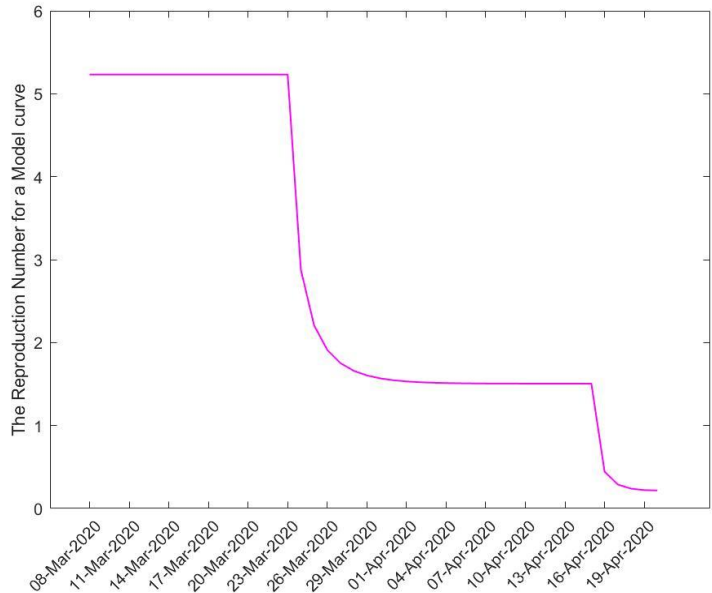

8. On Figure 1 we see that the nodes of the spline satisfy Node1 = **T2** = 17-Mar-2020 and Node2 = **T3** = 7-Apr-2020. On the other hand on Figure 2 we have chosen a configuration with different nodes, Node1 = **T2** = 23-Mar-2020, and Node2 = **T3** = 15-Apr-2020.
9. In the above examples of the dynamics of the Beta and Gamma coefficients; one sees the exponential factor $\exp(0.4 * (t - t_1))$ by which the curve changes from one level at $t = t_1$ to the next target level.
10. The functions $\beta(t)$ and $\gamma(t)$, are defined as interpolation splines on the subintervals defined by the Start date **T1**, the nodes dates **T2**, **T3**, and the final date **T4**. Thus the whole configuration is defined by eight parameters in total, which we gather in a set Θ, given by
    Θ = {T2, T3, β(T2), β(T3), β(T4), γ(T2), γ(T3), γ(T4)} .
11. The data which we use for the fitting of the TVBG-SEIR model are the official data for daily new infected cases *Idata(t)*, and the cumulative data for recovered and fatalities.
12. Finally, we fit the Model to the data by optimizing the positions of the two nodes **T2**, **T3**, and the levels of *β*(t) and *γ(t),* i.e. by applying Least squares by varying the set Θ in the usual quadratic functional $F(\Theta)$. Here, the nodes T2, T3 of the splines are varied from T1 to T4, and the interpolation values *β(T2), β(T3), β(T4), γ(T2), γ(T3), γ(T4),* are also independently varied. The model variables *I(t), R(t),* are obtained by solving the system (2a)-(2d) with initial conditions given by S(0) = N – I(0) – R(0), where N is the size of the whole population, *I(0) = Idata(0)*, and *R(0) = 0,* and the set Θ defines the splines for *β(t)* and *γ(t).* Additionally, one may introduce non-negative weights $w_1(t)$ and $w_2(t)$ which give priority to some of the data.
13. On the Figure below we provide the data *Idata(t)* for Bulgaria:

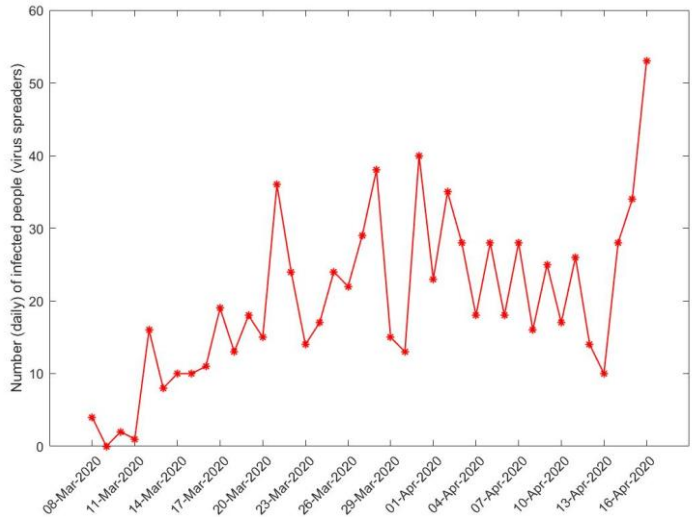

14. On the Figure below are provided the data for recovered and fatalities:

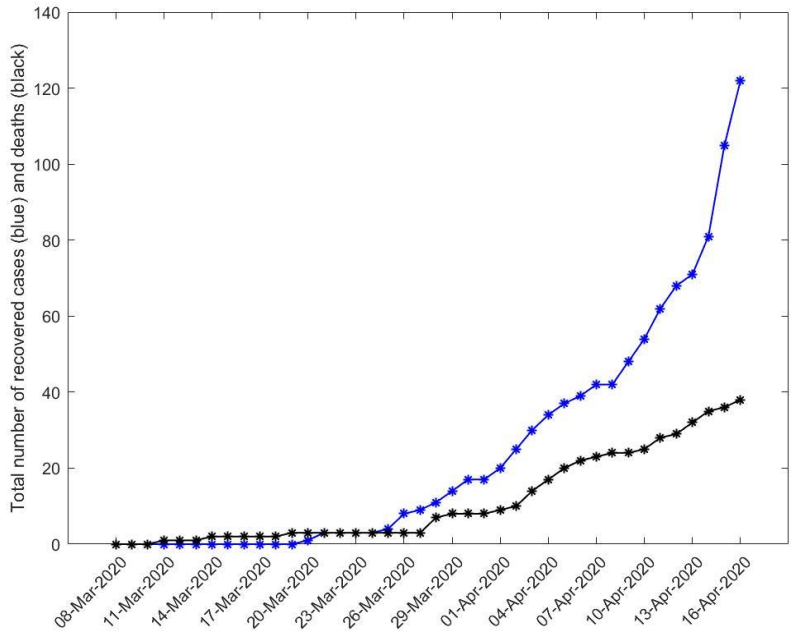

15. The result of the above minimization of the quadratic functional $F(\Theta)$ is a set of parameters $\bar{\Theta}$ for which the minimum is attained, whereby there may be multiple solutions.
On the Figure below we provide the fitting of the model curve *I(t)* to the above data for Bulgaria, *Idata(t)*:

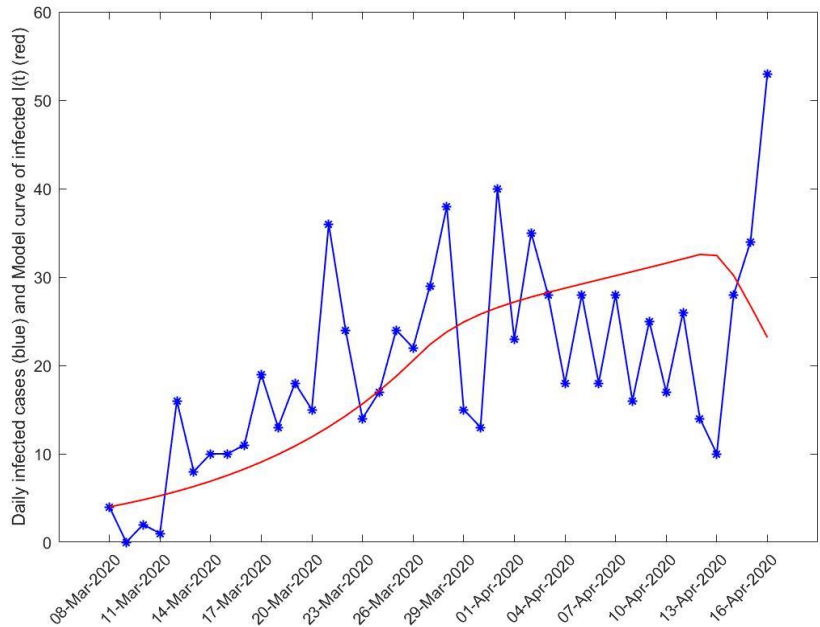

16. The Figure below shows the fitting by the model curve R(t) of the recovered and fatalities data for Bulgaria:

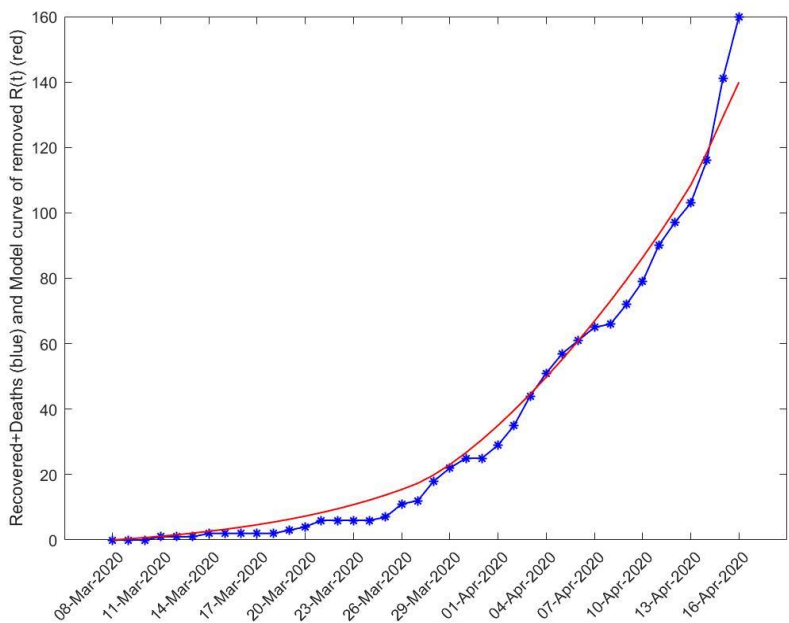

17. One of the most important tests for the quality of the model is to see how good does the optimized model identify the dates of introduction by the authorities of the Restrictive measures. This is clearly demonstrated below on the data for Bulgaria.
18. Final remark about the "parsimonious" style of constructing the spline model: one has to avoid putting too many nodes in the splines since this will influence the stability of the model, and might cause overfitting, hence would spoil the predictive power of the model.

## 4. Application of the model to Bulgarian COVID-19 data and scenarios from the TOOL

Here we demonstrate how the Tool works in the case of Bulgarian data. The Tool will be described in detail in Section 7.

We provide the visualizations of the model fitting which are available in the Tool. The thick red line shows (until Today = T4) the fitted model curve for the daily new infected cases *I(t)* and the blue stars are the official data *Idata(t).* The thin red line shows the prediction scenarios, after Today.

1. For Bulgaria we have considered the data from the Start Date which is **T1** = 8-Mar-2020, until the end date Today = **T4** = 18-Apr-2020. The minimization of the functional $F(\Theta)$ has resulted in *top three* Models: **Model1**, **Model2** and **Model3**, defined by the sets $\bar{\Theta}$ $\overline{\Theta^0}$, $\overline{\Theta^i}$, for which the value ***Fval*** of the functional $F(\Theta)$ is very close to the minimum (which is attained by the set $\bar{\Theta}$ for Model1). Let us note that all $F(\Theta)$ have the **maximum** value $0.1262*10^3$ (taken over all admissible parameters $\Theta$).

2. It is an important result that the first node of **Model1** is on 18-Mar-2020 which is between the First Restr. Measures date 14-Mar-2020 and the Second Restr. Measures date 20-Mar-2020.
   This finding shows that our model reflects the realities.

3. As is seen from the Figures below, the top three models indicate the peak of the curve which is in the vicinity of the time interval 12-Apr-2020 and 15-Apr-2020.

4. For **Model1** we have $F(\bar{\Theta}) = 0.1043*10^3$, hence the ratio max/$F(\bar{\Theta})$ is about 1.21. The curve *I(t)* of Model1 shows a peak at the nodal point 12-Apr-2020.
   The Figure below shows the **simplest prediction scenario** starting on Today = 18-Apr-2020. In the Legend of the Figure, Coef1 = 1 and Coef2 = 1 mean that no change by the authorities will be undertaken starting Today and ending on 13-May-2020. Further, Coef11 = 1 and Coef22 = 1 mean that no relaxation of the measures will follow starting on 13-May-2020.

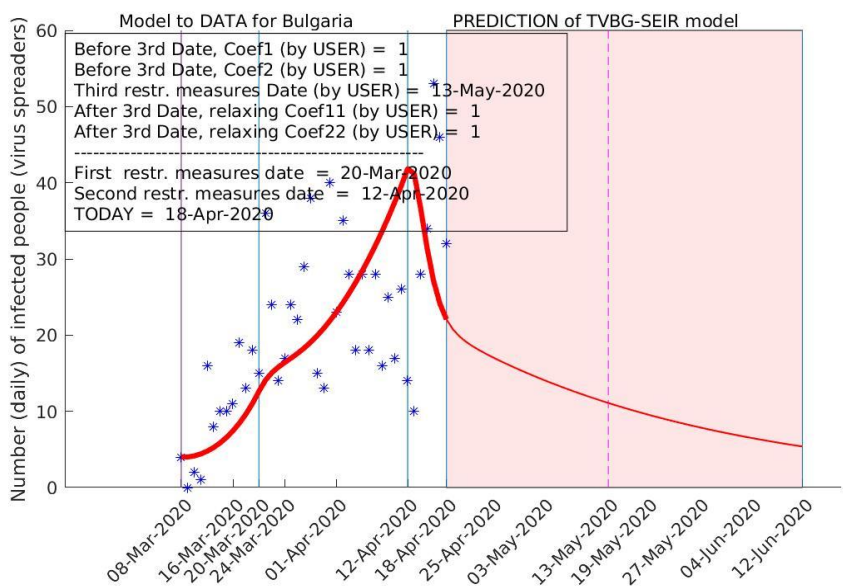

5. However, on 13-May-2020 only the Beta measures may be partially relaxed (Coef11 = 1.4), bringing to 10 infected cases per day. This is seen on the Figure below:

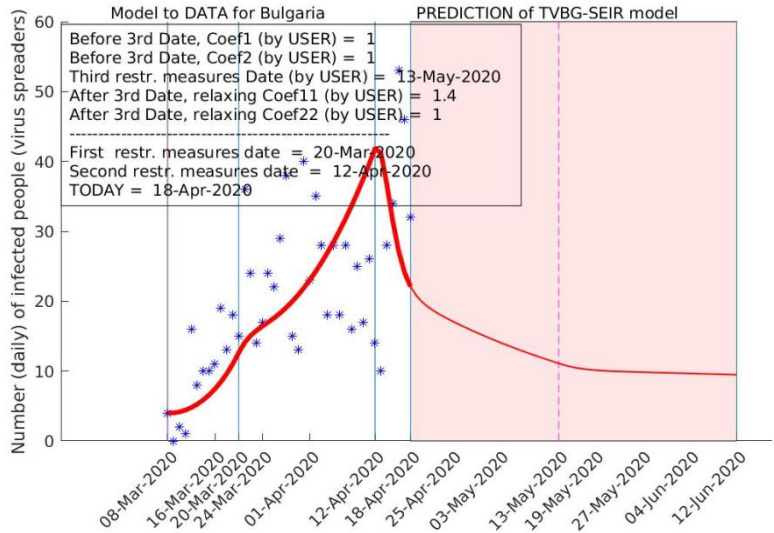

6. But these measures cannot be relaxed more, e.g. with Coef11 = 1.8, as seen from the Figure below:

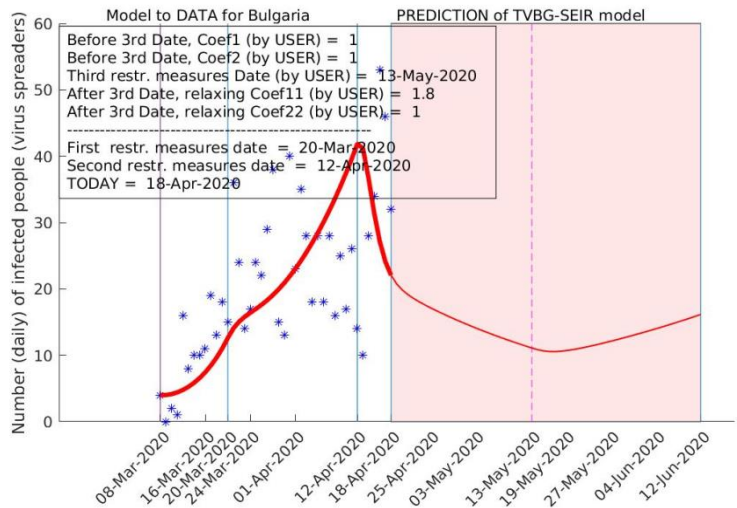

7. Also, relaxing both measures (Coef11=Coef22=1.4) after 13-May-2020 will not be good, as seen from the Figure below:

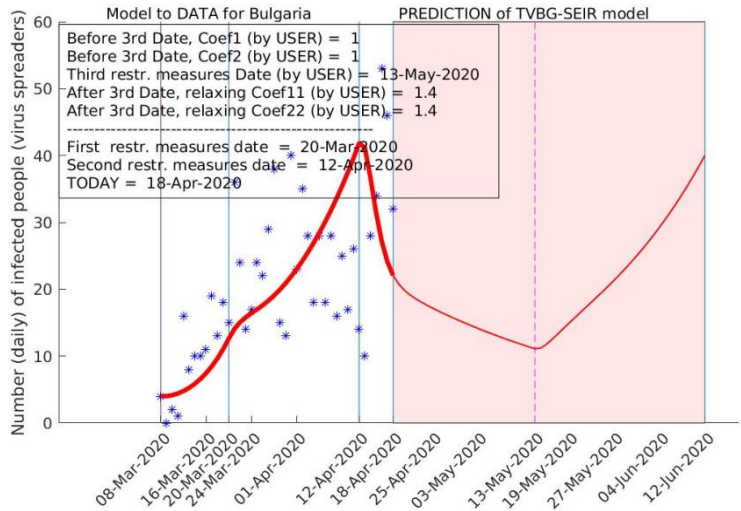

Similar are the conclusions with **Model3** for Bulgaria (with Fval=0.1070*103 )

8. For **Model2** (with Fval=0.1057*103 ) we have the most optimistic scenario since we may partially relax both measures after 13-May-2020 (i.e. Coef11=Coef22=1.4):

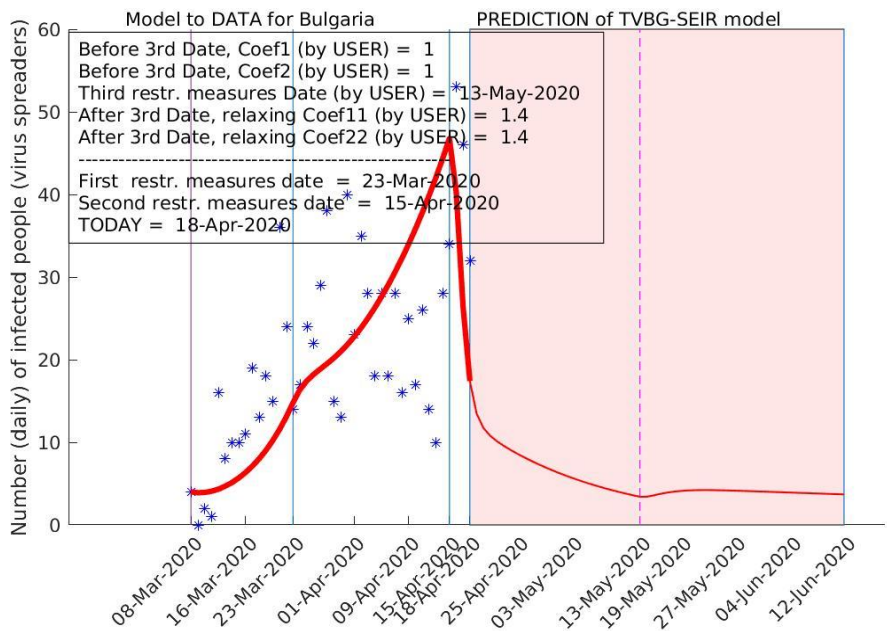

Obviously, from practical point of view one has to take into account the predictions of the three models.

## 5. Application of the model to Italian COVID-19 data and scenarios from the TOOL

1. The results about **Italy** considered till Today = 18-Apr-2020 are quite unusual.
Although the official data show that there is a serious decline, our curves are reluctant to go down. Apparently, the simultaneously fitted Removed data (containing Recovered cases and Fatalities) are pulling the curve up. None of the Models which are suboptimal, has nodes which are in April. This is very untypical.

2. We have to go into details of our method to explain things better. The maximum value of the Quadratic Functional (which is minimized in the LS fitting) is 103 *113.015 and the ration max/min is about 10. The value of the functional for the first 5 best Models (close to the minimal value) are: 10.527, 10.579, 10.627, 10.717, 10.977, all multiplied by 103, which shows that they give also a good approximation of the real data. We showcase the first three models where no control of the measures for scenarios is considered (i.e. Coef1=Coef2=Coef11=Coef22=1):

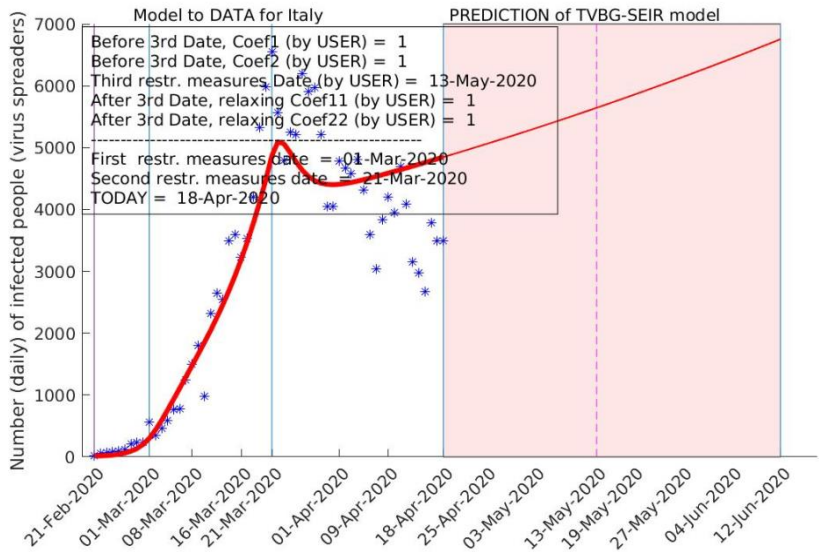

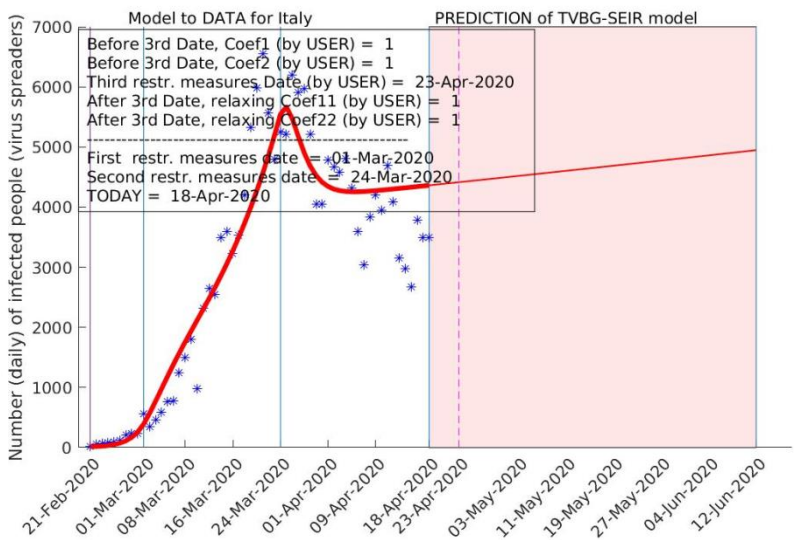

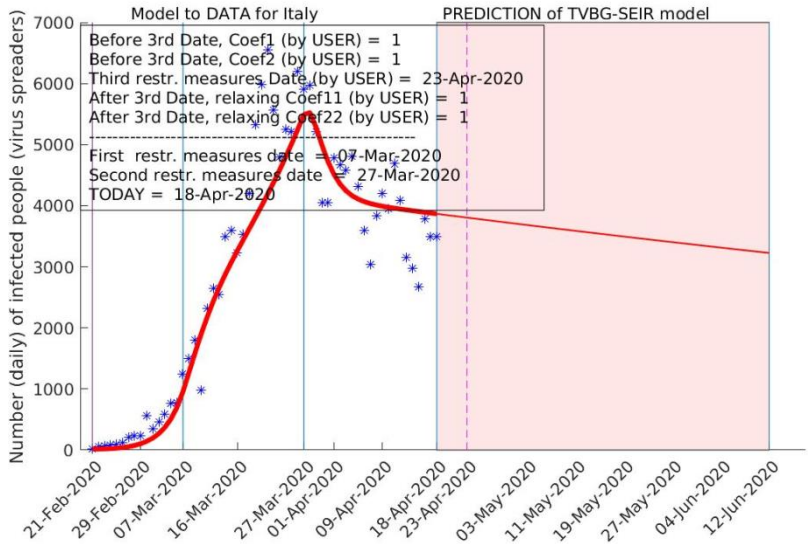

3. We decided to make a special test for the data in Italy. Since the number of fatalities is very essential, and it is notorious that the number of unreported cases may be even 100 times more than the reported, it is worth doing the following experiment – just multiply the Daily infected by 100, and the Recovered by 100, and then fit the Model and see the result. The surprise is that it gave different result but still the reluctance of the curve I(t) to go down has remained:
**Model1** is the most optimistic:

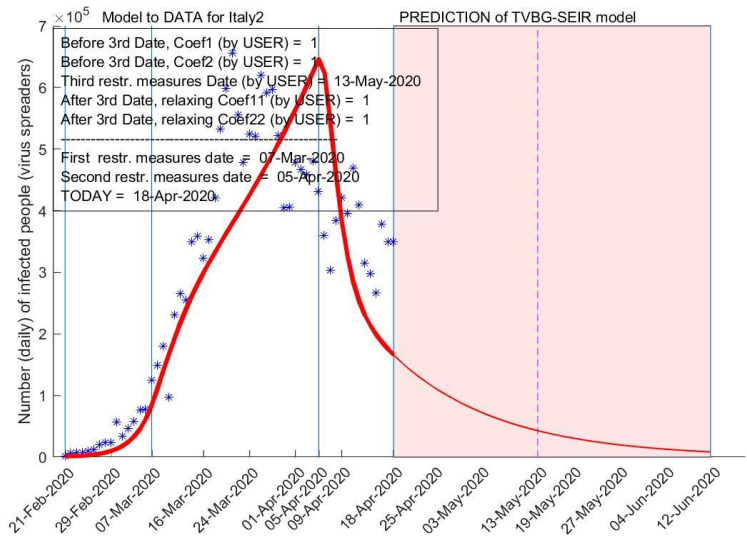

**Model2** is not that optimistic:

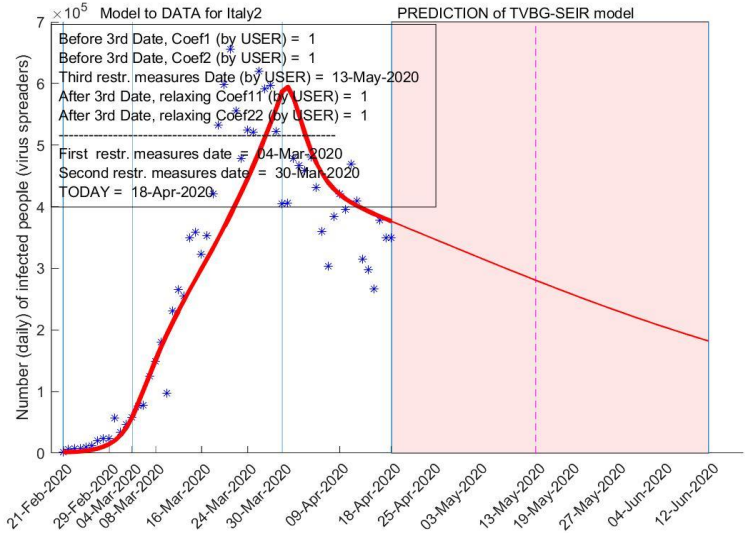

And **Model3** is also not that optimistic:

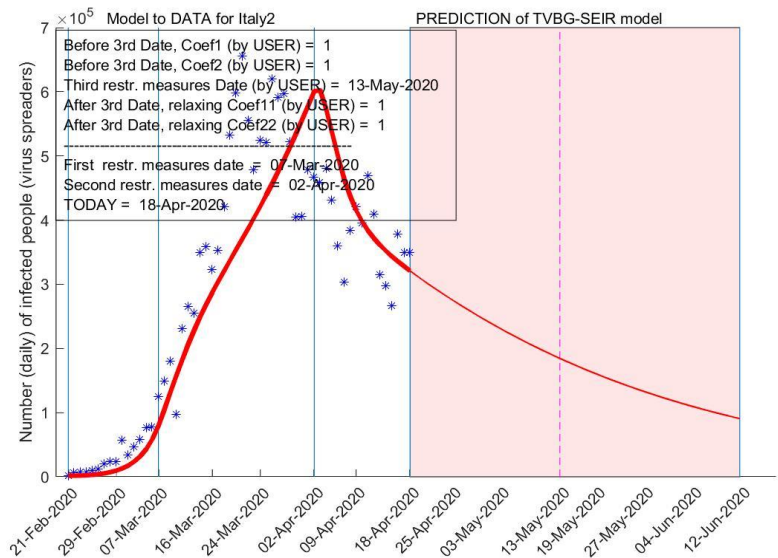

## 6. Application of the model to German COVID-19 data and scenarios from the TOOL

1. For the **German** data, on the same Today date = 18-Apr-2020, according to our **Model1** for Germany (with Fval=20.849*103 , where the maximum value of Fval is 118.482*103, hence the ratio max/min is about 6), one may relax safely just a little bit the Beta measures (Coef11 = 1.4), bringing to less than 1000 infected per month:

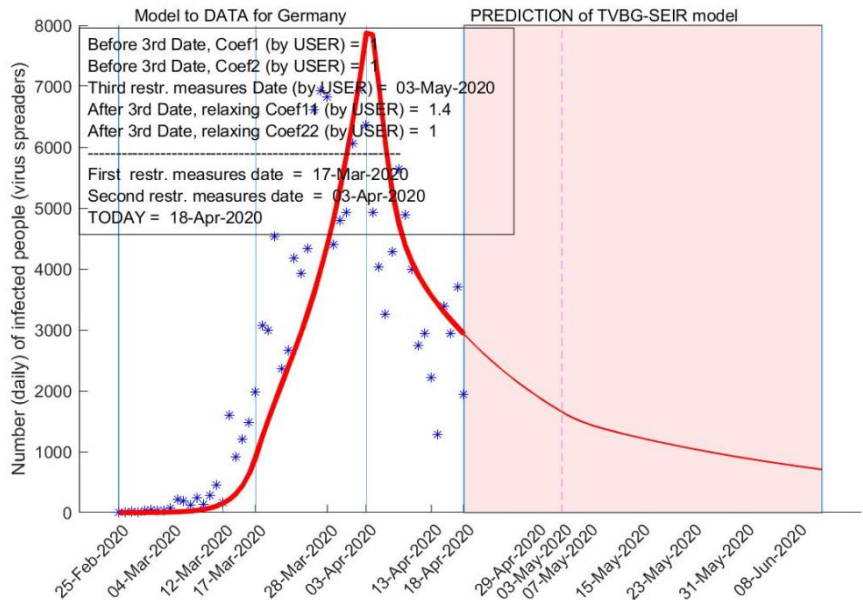

2. however not both measures could be relaxed (Coef11=Coef22=1.4):

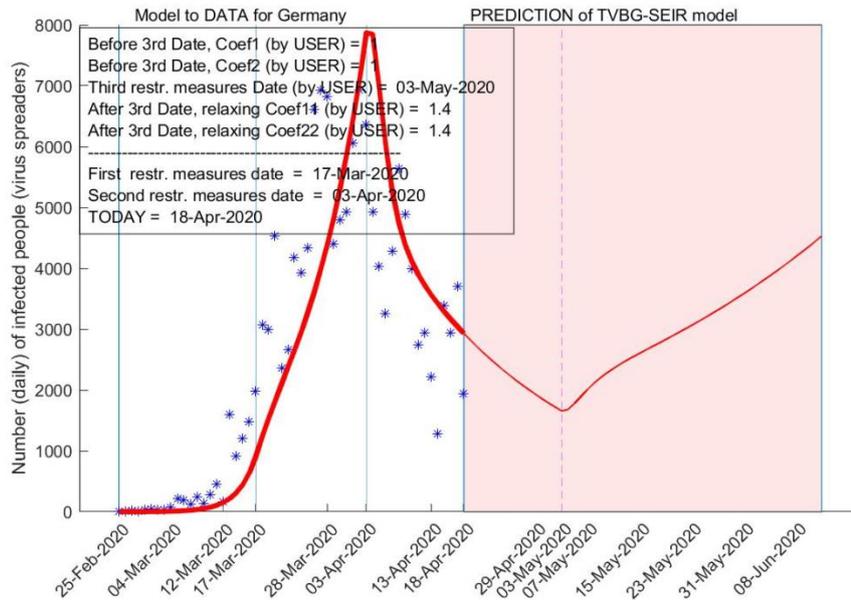

3. On the other hand, **Model2** (which is somewhat less reliable since it fits a little worse the data, Fval=21.368*103) with both measures relaxed (Coef11=Coef22=1.4) gives again less than 1000 infected per day:

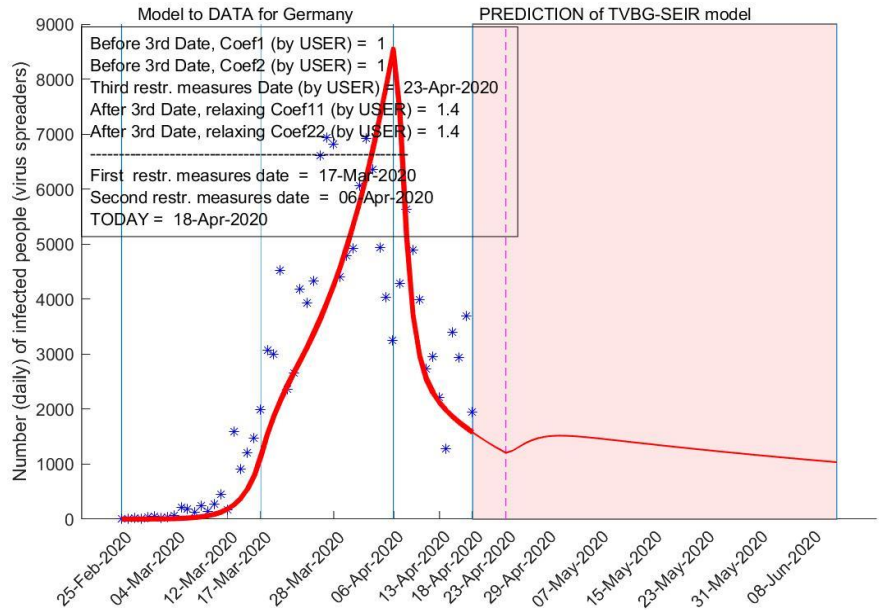

## 7. Description of the TOOL for controlled scenarios

1. We have designed a TOOL for the Model predictions of the Coronavirus (and similar infectious diseases) spread. The tool is based on the fitting of the TVBG-SEIR model to the official data described above.
   It is available at the site
   http://213.191.194.141:8888/notebooks/TVBG-SEIR-Spline-model_v3.ipynb?token=b5d97bfbd7dd062e47ee7ab51837e470a8c226743a4667ee
2. First of all, on Today's date **T4**, we fit the model for the time series in the interval [**T1**,**T4**]. Then the USER may choose several parameters to make a prediction about the virus spread during the period [T4, Horizon] where the Horizon is chosen to be about 2 months from Today (=T4).
3. The first parameter, to be controlled, is called Third restrictive measures date denoted by **T5**. The USER may choose several options, say 5, 15, 25 days from Today (=**T4**), i.e. one may select the dates **T5 =** T4+5, T4+15, T4+25.
4. Then the USER may decide how to strengthen or relax the Beta measures and the Gamma measures during the period [**T4**,**T5**], by means of the coefficients Coef1 and Coef2 respectively; Coef1 = 1 means that the Beta measures remain the same, while Coef2 = 1 means that the Gamma measures remain the same in the period [**T4**,**T5**]. If Coef1 is less than 1 then this means that the Beta measures are "weaker", and also, the smaller Coef1, the weaker are the Beta measures and they will reach a target value at the date **T5**, which is defined by the size of Coef1 (Note that Coef1 < 1 means that the Beta rate will be bigger!). In a similar way, if Coef2 is less than 1, then this means that the Gamma measures will be "weaker", and the smaller Coef2, the weaker are the Gamma measures (Note that then the Gamma rate will be smaller!). A target value (determined by the size of Coef2) will be reached at the date **T5**.
5. On the other hand, if Coef1 or Coef2 are bigger than 1, this means "strengthening the measures", resp. of Beta measures and Gamma measures in the period [**T4**,**T5**] to some target value defined by Coef1, Coef2.

6.  The USER has further the possibility to decide what will happen after date **T5** - to weaken (or leave the same) the Beta and the Gamma measures. This is decided by the choice of two coefficients – Coef11 for the Beta and Coef22 for the Gamma measures. Coef11 = 1 means that one retains the same level of the Beta measures; Coef22 = 1 means that one retains the same level of the Gamma measures. If Coef11 is bigger than 1 then this would relax the Beta measures – the bigger Coef11 the more the relaxation. Coef22 makes the same for the Gamma measures.

## 8. Recent research on time-varying transmission rates

As we already said, presently it is urgent to consider SIR/SEIR models with time-varying $\beta(t)$ and $\gamma(t)$ rates.

Let us mention some research about solving an inverse problem for finding time-varying $\beta(t)$ in a SIR model, for a fixed removal rate $\gamma$, [10], where the time-varying transmission rate $\beta(t)$ is determined by the infectious cases. In [1], the authors do research and provide further references of research on specific models for the transmission rate $\beta(t)$.

## Acknowledgements

All authors acknowledge the support by grants DH-02-13 and KP-06-N32-8 of Bulgarian NSF. The first-named author acknowledges the partial support by Grant No BG05M2OP001-1.001-0003, financed by the Science and Education for Smart Growth Operational Program (2014-2020) and co-financed by the European Union through the European structural and Investment funds.